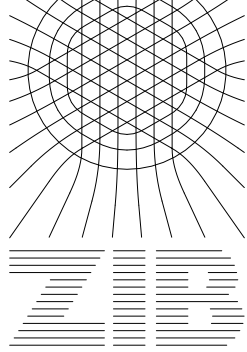

Wolfram Koepf

# A Package on Orthogonal Polynomials and Special Functions

# A Package on Orthogonal Polynomials and Special Functions

Wolfram Koepf

Konrad-Zuse-Zentrum für Informationstechnik

Takustr. 7

D-14195 Berlin

email: `koepf@zib.de`

**Abstract:**

In many applications (hypergeometric-type) special functions like orthogonal polynomials are needed. For example in more than 50 % of the published solutions for the (application-oriented) questions in the "Problems Section" of SIAM Review special functions occur.

In this article the Mathematica package `SpecialFunctions` which can be obtained from the URL `http://www.zib.de/koepf` is introduced [15]. Algorithms to convert between power series representations and their generating functions is the main topic of this package ([8]–[15]), extending the previous package `PowerSeries` [12]. Moreover the package automatically finds differential and recurrence equations ([13]–[14]) for expressions and for sums (the latter using Zeilberger's algorithm ([23], [18], [13]).

As an application the fast computation of polynomial approximations of solutions of linear differential equations with polynomial coefficients is presented. This is the asymptotically fastest known algorithm for series computations, and it is much faster than Mathematica's builtin `Series` command if applicable. Many more applications are considered.

Finally the package includes implementations supporting the efficient computation of classical continuous and discrete orthogonal polynomials.

## 1 Holonomic Functions

A homogeneous linear differential equation

$$\sum_{k=0}^{m} p_k(x) f^{(k)}(x) = 0$$

with polynomial coefficients $p_k(x)$ is called *holonomic*, as is the corresponding $f(x)$.

Holonomic functions have nice algebraic properties: In 1980 Stanley [22] proved by algebraic arguments that the sum and product of holonomic functions, and their composition with algebraic functions, in particular with rational functions and rational powers, form holonomic functions, again.

By iterative differentiation and the use of Gaussian elimination, one can construct the resulting differential equations algorithmically.

The funny thing is that these algorithms had been known already in the last century ([1]–[2]), but because of their complexity had fallen into oblivion.

The package `SpecialFunctions` contains implementations of these algorithms. The Mathematica procedure `SumDE[DE1,DE2,f[x]]` computes the holonomic differential equation for the sum of two functions which correspond to the holonomic differential equations `DE1` and `DE2`, written in terms of `f[x]`. Here are some examples:

```
In[1]:= <<SpecialFunctions`
SpecialFunctions, (C) Wolfram Koepf, version 1.00, November 20, 1996
Fast Zeilberger, (C) Peter Paule and Markus Schorn (V 2.2) loaded

In[2]:= SumDE[f''[x]+f[x]==0,f'[x]-f[x]==0,f[x]]

                                   (3)
Out[2]= f[x] - f'[x] + f''[x] - f   [x] == 0

In[3]:= ProductDE[f''[x]+f[x]==0,f'[x]-f[x]==0,f[x]]

Out[3]= 2 f[x] - 2 f'[x] + f''[x] == 0

In[4]:= SumDE[f''[x]-x*f[x]==0,f'[x]-f[x]==0,f[x]]

                 2
Out[4]= (1 - x + x ) f[x] + (1 - x) x f'[x] -

                            (3)
>       x f''[x] + (-1 + x) f   [x] == 0

In[5]:= ProductDE[f''[x]-x*f[x]==0,f'[x]-f[x]==0,f[x]]

Out[5]= (-1 + x) f[x] + 2 f'[x] - f''[x] == 0
```

In every day language the first computation reads as: The sum of the functions $\sin x$ (satisfying $f''(x)+f(x)=0$) and $e^x$ (satisfying $f'(x)-f(x)=0$), $\sin x + e^x$, satisfies the differential equation $f(x)-f'(x)+f''(x)-f'''(x)=0$. Actually any linear combination $a\sin x + b\cos x + ce^x$ satisfies this resulting differential equation. In the next line the differential equation for the product $\sin x\, e^x$ (or $\cos x\, e^x$) is computed. Note that the further computations are rather similar, although one of the corresponding functions is much more advanced, namely an Airy function (satisfying $f''(x)-x f(x)=0$).

If functions are given by expressions, holonomic differential equations can be determined as well. The procedure `HolonomicDE[expr,f[x]]` internally uses the sum and product algorithms:

```
In[6]:= HolonomicDE[AiryAi[x]*Exp[x],f[x]]

Out[6]= (-1 + x) f[x] + 2 f'[x] - f''[x] == 0

In[7]:= HolonomicDE[Exp[alpha x]*Sin[beta x],f[x]]

             2      2
Out[7]= (alpha  + beta ) f[x] - 2 alpha f'[x] +

>       f''[x] == 0

In[8]:= de=HolonomicDE[ArcSin[x],f[x]]
                          2
Out[8]= x f'[x] + (-1 + x ) f''[x] == 0

In[9]:= Nest[ProductDE[de,#,f[x]]&,de,4]
```

```
Out[9]= x f'[x] + (-16 + 31 x ) f''[x] +

                    2    (3)
>      15 x (-5 + 6 x ) f   [x] +

                             2   (4)
>      5 (1 - x) (1 + x) (4 - 13 x ) f   [x] +

                2        2  (5)
>      15 (-1 + x)  x (1 + x)  f   [x] +

              3        3  (6)
>      (-1 + x)  (1 + x)  f   [x] == 0
```

The last statement computes the holonomic differential equation for $\arcsin^5 x$. This result can also be obtained by the simple statement `HolonomicDE[ArcSin[x]^5,f[x]]`.

## 2 Series Representations

By an easy procedure (equating coefficients) any holonomic differential equation can be converted to a holonomic recurrence equation for a series representation

$$f(x) = \sum_{k=-\infty}^{\infty} a_k \, x^k \; .$$

It turns out that a function $f(x)$ is holonomic if and only if its series coefficients $a_k$ are holonomic, i.e., if the latter satisfy a homogeneous, linear recurrence equation

$$\sum_{k=0}^{M} a_k(n) \, f_{n+k}(x) = 0$$

with polynomial coefficients $a_k$ in $n$. The procedure `DEtoRE[DE,f[x],a[k]]` converts the differential equation `DE` for $f(x)$ to the corresponding holonomic recurrence equation for the coefficients $a_k$, and `REtoDE[RE,a[k],f[x]]` does the corresponding back conversion. We have for example

```
In[10]:= de=HolonomicDE[ArcSin[x]^2,f[x]]

Out[10]= f'[x] + 3 x f''[x] +

                        (3)
>      (-1 + x) (1 + x) f   [x] == 0

In[11]:= re=DEtoRE[de,f[x],a[k]]

          3
Out[11]= k  a[k] - k (1 + k) (2 + k) a[2 + k] == 0

In[12]:= REtoDE[re,a[k],f[x]]
```

```
Out[12]= f'[x] + 3 x f''[x] - f   [x] + x  f   [x] == 0
                         (3)     2  (3)
```

If for a given input function $f(x)$ the resulting recurrence equation has only two summands, one has a *generalized hypergeometric function* (see e.g. [21]), and the coefficients can be computed explicitly.

The generalized hypergeometric series is given by

$$_pF_q\left(\begin{array}{cccc} a_1 & a_2 & \cdots & a_p \\ b_1 & b_2 & \cdots & b_q \end{array} \bigg| x\right) := \sum_{k=0}^{\infty} A_k\, x^k$$

$$= \sum_{k=0}^{\infty} \frac{(a_1)_k \cdot (a_2)_k \cdots (a_p)_k}{(b_1)_k \cdot (b_2)_k \cdots (b_q)_k \, k!} x^k$$

where $(a)_k := \prod_{j=1}^{k}(a{+}j{-}1) = \Gamma(a{+}k)/\Gamma(a)$ denotes the *Pochhammer symbol* (`Pochhammer[a,k]`) or *shifted factorial*. $A_k$ is a *hypergeometric term* and fulfils the recurrence equation ($k = 0, 1, \ldots$)

$$A_{k+1} := \frac{(k+a_1)\cdots(k+a_p)}{(k+b_1)\cdots(k+b_q)(k+1)} \cdot A_k$$

with the initial value

$$A_0 := 1\,.$$

In Mathematica the generalized hypergeometric function is denoted by
`HypergeometricPFQ[plist,qlist,x]`, where

$$\text{plist} = [a_1, a_2, \ldots, a_p] \quad \text{and} \quad \text{qlist} = [b_1, b_2, \ldots, b_q]\,.$$

To find the representing Laurent-Puiseux expansion of a given expression $f$ with respect to the variable $x$ and point of development $x_0$ one uses the function call `PowerSeries[f,{x,x0}]` (short `PS[f,{x,x0}]`). The package also supports the syntax `PowerSeries[f,x,x0]`. If $x_0 = 0$ one can also use either `PowerSeries[f,x]`, or `PS[f,x]`. Furthermore the procedure `FunctionToHypergeometric[f,x]` brings the resulting series in hypergeometric notation. One has e.g.

```
In[13]:= PS[ArcSin[x]^2,x]

              k  2 + 2 k      2
             4  x         k!
Out[13]= sum[------------------, {k, 0, Infinity}]
              (1 + k) (1 + 2 k)!

In[14]:= FunctionToHypergeometric[ArcSin[x],x]

                      1   1     3     2
Out[14]= x hypergeometricPFQ[{-, -}, {-}, x ]
                      2   2     2

In[15]:= FunctionToHypergeometric[ArcSin[x]^2,x]
```

```
                 2                                          3    2
Out[15]= x  hypergeometricPFQ[{1, 1, 1}, {2, -}, x ]
                                            2

In[16]:= PS[LaguerreL[n,x],x]

              k
             x  Pochhammer[-n, k]
Out[16]= sum[-------------------, {k, 0, Infinity}]
                       2
                      k!
```

Here `LaguerreL[n,x]` denotes the $n^{th}$ Laguerre polynomial $L_n(x)$.

We use the head `sum` to denote a resulting infinite sum since Mathematica Version 3 evaluates sums with head `Sum` instantly. Similarly Mathematica's builtin procedure `HypergeometricPFQ` has instant evaluations in some instances, e.g.

```
In[17]:= HypergeometricPFQ[{1,1},{2},x]

           Log[1 - x]
Out[17]= -(----------)
               x
```

By this reason we use the (non-evaluating) name `hypergeometricPFQ` instead.

Note that the procedure `PowerSeries` can handle the more general case of *Laurent-Puiseux* series representations $(p = 1, 2, \ldots)$

$$f(x) = \sum_{k=-\infty}^{\infty} a_k\, x^{k/p},$$

for example

```
In[18]:= PS[Sin[Sqrt[x]],x]

               k  1/2 + k
           (-1)  x
Out[18]= sum[--------------, {k, 0, Infinity}]
              (1 + 2 k)!
```

## 3  Recurrence Equations

Similarly to the sum and product algorithms for holonomic differential equations there are simple sum and product algorithms for holonomic recurrence equations. These compute the recurrence equation for the sum and product, respectively, of two discrete functions $a_k$ and $b_k$, given by two holonomic recurrence equations `RE1` and `RE2`. These algorithms are invoked by the procedures `SumRE[RE1,RE2,a[k]]` and `ProductRE[RE1,RE2,a[k]]`, respectively. We get for example

```
In[19]:= SumRE[a[k+1]==a[k]/(k+1),a[k+1]==a[k],a[k]]
                                  2
Out[19]= (1 + k) a[k] + (-1 - 3 k - k ) a[1 + k] +

>    k (2 + k) a[2 + k] == 0

In[20]:= ProductRE[a[k+1]==a[k]/(k+1),a[k+1]==a[k],a[k]]

Out[20]= a[k] + (-1 - k) a[1 + k] == 0
```

Here the first computation finds the recurrence equation valid for the sum of $1/k!$ (satifying $a_{k+1} = a_k/(k+1)$) and 1 (satifying $a_{k+1} = a_k$), and the second one for the corresponding product. From this point of view the result of the second computation is obvious.

The automatic computation of holonomic recurrence equations for expressions is supported through the procedure `HolonomicRE[expr,a[k]]`. Examples are

```
In[21]:= HolonomicRE[(n!+k!^2)/k,a[n]]

              2                    2
Out[21]= (1 + n)  a[n] + (-1 - 3 n - n ) a[1 + n] +

>    n a[2 + n] == 0

In[22]:= HolonomicRE[(n!+k!^2)/k,a[k]]

              3
Out[22]= k (1 + k)  (3 + k) a[k] -

                  2              2
>    (1 + k) (1 + 3 k + k ) (3 + 3 k + k ) a[1 + k] +

              2
>    k (2 + k)  a[2 + k] == 0
```

The work with special functions for which holonomic recurrence equations exist, is possible.

```
In[23]:= HolonomicRE[LegendreP[k,x],P[k]]

Out[23]= (1 + k) P[k] - (3 + 2 k) x P[1 + k] +

>    (2 + k) P[2 + k] == 0

In[24]:= HolonomicRE[Binomial[2*k,k]*HermiteH[k,x],a[k]]

Out[24]= 8 (1 + 2 k) (3 + 2 k) a[k]

>    - 4 (3 + 2 k) x a[1 + k] + (2 + k) a[2 + k] == 0
```

One has also access to *Zeilberger's algorithm*. For this purpose the package uses an implementation of Paule/Schorn from RISC [18]. Zeilberger's algorithm finds a holonomic recurrence equation for a sum

$$S_n := \sum_{k=k_1}^{k_2} F(n,k)$$

if $F(n,k)$ is a hypergeometric term w.r.t. both $n$ and $k$, i.e., the term ratios

$$\frac{F(n+1,k)}{F(n,k)} \in \mathbb{Q}(n,k), \qquad \frac{F(n,k+1)}{F(n,k)} \in \mathbb{Q}(n,k)$$

are both rational functions in $n$ and $k$. One gets for example for the binomial powers

```
In[25]:= HolonomicRE[sum[Binomial[n,k]^5,{k,0,n}],a[n]]

                    4                   2
Out[25]= 32 (1 + n)  (292 + 253 n + 55 n ) a[n] +

                                              2              3
>        (-514048 - 1827064 n - 2682770 n  - 2082073 n  -

                 4            5            6
>        900543 n  - 205799 n  - 19415 n ) a[1 + n] +

                                             2              3
>        (-79320 - 245586 n - 310827 n  - 205949 n  -

                 4           5           6
>        75498 n  - 14553 n  - 1155 n ) a[2 + n] +

                 4                   2
>        (3 + n)  (94 + 143 n + 55 n ) a[3 + n] == 0
```

Some time ago, such a result was worth a publication ([19], [4]).
Similarly, for the *Apéry numbers*

$$A_n := \sum_{k=0}^{n} \binom{n}{k}^2 \binom{n+k}{n}^2 \tag{1}$$

one deduces the recurrence equation

```
In[26]:= HolonomicRE[sum[Binomial[n,k]^2*
            Binomial[n+k,k]^2,{k,0,n}],A[n]]

                3
Out[26]= (1 + n)  A[n] -

                                     2
>        (3 + 2 n) (39 + 51 n + 17 n ) A[1 + n] +

                3
>        (2 + n)  A[2 + n] == 0
```

which played an essential role in Apéry's proof (see e.g. van der Poorten's [20] entertaining presentation of Apéry's proof and its history[1]) of the irrationality of

$$\zeta(3) = \sum_{k=1}^{\infty} \frac{1}{k^3}.$$

---

[1] Van der Poorten writes: "To convince ourselves of the validity of Apéry's proof we need only complete the following exercise: show that the recurrence equation is valid for (1). Neither Cohen nor I had been able to prove this in the intervening two months..."

Zeilberger's algorithm can be used to prove identities. The *Legendre polynomials* $P_n(x)$, e.g., can be represented as any of

$$\begin{aligned} P_n(x) &= \sum_{k=0}^n \binom{n}{k} \binom{-n-1}{k} \left(\frac{1-x}{2}\right)^k \\ &= {}_2F_1\left(\begin{array}{c} -n, n+1 \\ 1 \end{array} \bigg| \frac{1-x}{2}\right) \\ &= \frac{1}{2^n} \sum_{k=0}^n \binom{n}{k}^2 (x-1)^{n-k} (x+1)^k \\ &= \left(\frac{1-x}{2}\right)^n {}_2F_1\left(\begin{array}{c} -n, -n \\ 1 \end{array} \bigg| \frac{1+x}{1-x}\right) \\ &= \frac{1}{2^n} \sum_{k=0}^{\lfloor n/2 \rfloor} (-1)^k \binom{n}{k} \binom{2n-2k}{n} x^{n-2k} \\ &= \binom{2n}{n} \left(\frac{x}{2}\right)^n {}_2F_1\left(\begin{array}{c} -n/2, -n/2+1/2 \\ -n+1/2 \end{array} \bigg| \frac{1}{x^2}\right) \\ &= x^n \, {}_2F_1\left(\begin{array}{c} -n/2, -n/2+1/2 \\ 1 \end{array} \bigg| 1 - \frac{1}{x^2}\right). \end{aligned}$$

By the following computations three of these representations are proved to be consistent

```
In[27]:= HolonomicRE[sum[Binomial[n,k]*
         Binomial[-n-1,k]*((1-x)/2)^k,{k,0,n}],P[n]]

Out[27]= (-1 - n) P[n] + (3 + 2 n) x P[1 + n] +

>        (-2 - n) P[2 + n] == 0

In[28]:= HolonomicRE[1/2^n*sum[Binomial[n,k]^2*
         (x-1)^(n-k)*(x+1)^k,{k,0,n}],P[n]]

Out[28]= (1 + n) P[n] - (3 + 2 n) x P[1 + n] +

>        (2 + n) P[2 + n] == 0

In[29]:= HolonomicRE[1/2^n*sum[(-1)^k*Binomial[n,k]*
         Binomial[2n-2k,n]*x^(n-2k),{k,0,Infinity}],P[n]]

Out[29]= (-1 - n) P[n] + (3 + 2 n) x P[1 + n] +

>        (-2 - n) P[2 + n] == 0
```

To finish the proof that these sums represent the same family of functions, it is enough to verify that two initial values agree since by the above computations the sums satisfy the same second order recurrence equation. Note that the last sum is supported in the interval $[0, \lfloor n/2 \rfloor]$ which can be seen from the two upper indices $-\frac{n-1}{2}$ and $-\frac{n}{2}$ of the hypergeometric representation

```
In[30]:= SumToHypergeometric[sum[1/2^n*(-1)^k*Binomial[n,k]*
         Binomial[2n-2k,n]*x^(n-2k),{k,0,Infinity}]]

            n
Out[30]= (x  Gamma[1 + 2 n]

                         1 - n   -n    1 - 2 n      -2
>        hypergeometricPFQ[{-----, --}, {-------}, x  ])\
                           2     2        2

           n           2
>        / (2  Gamma[1 + n]  )
```

The procedure `SumToHypergeometric[sum[expr,{k,0,Infinity}]]` converts an infinite sum into hypergeometric notation; the procedure `FunctionToHypergeometric` which was mentioned before uses `SumToHypergeometric` internally.

## 4  Generating functions

Generating functions can be calculated by converting the holonomic recurrence equation of the coefficient sequence to the corresponding differential equation for its generating function if applicable.

The procedure `Convert[sum[expr,{k,k0,Infinity}],x]` converts the Laurent-Puiseux series `Sum[expr,{k,k0,Infinity}]` w.r.t. $x$ to an expression if the corresponding differential equation can be solved by `DSolve`. For example, we get

```
In[31]:= Convert[sum[Binomial[n,k] x^k,{k,0,Infinity}],x]

                n
Out[31]= (1 + x)

In[32]:= Convert[sum[x^(2k+1)/(2k+1)!,{k,0,Infinity}],x]

               x
         -1    E
Out[32]= ---- + --
           x    2
         2 E
```

The procedure `GeneratingFunction[a,k,z]` tries to find the generating function

$$\sum_{k=0}^{\infty} a_k \, z^k \;,$$

and the procedure `ExponentialGeneratingFunction[a,k,z]` searches for the exponential generating function

$$\sum_{k=0}^{\infty} \frac{a_k}{k!} z^k \;.$$

Both functions use `Convert` internally. We get e.g. for the Legendre polynomials

```
In[33]:= specfunprint

In[34]:= GeneratingFunction[LegendreP[n,x],n,z]
SpecialFunctions, (C) Wolfram Koepf, version 1.00,
>    November 20, 1996
specfun-info: RE:
         (1 + k)*a[k] - (3 + 2*k)*x*a[1 + k] +
         (2 + k)*a[2 + k] == 0
specfun-info: DE:
                                              2
         (-x + z) F[z] + (1 - 2 x z + z ) F'[z] == 0
specfun-info: Trying to solve DE ...
specfun-info: DSolve computes
         C[1]/(1 - 2*x*z + z^2)^(1/2)
specfun-info: expression rearranged:
         C[1]/(1 - 2*x*z + z^2)^(1/2)
specfun-info: Calculation of initial values...
specfun-info: C[1] == 1

                  1
Out[34]= --------------------
                       2
         Sqrt[1 - 2 x z + z ]
```

The command `specfunprint` turns on a verbatim mode with which you will be informed about intermediate computation steps.

Similarly, we get

```
In[35]:= ExponentialGeneratingFunction[HermiteH[n,x],n,z]
SpecialFunctions, (C) Wolfram Koepf, version 1.00,
>    November 20, 1996
specfun-info: REProduct entered
specfun-info: RE:
         2*a[k] - 2*x*a[1 + k] + (2 + k)*a[2 + k] == 0
specfun-info: DE:
         2 (-x + z) f[z] + f'[z] == 0
specfun-info: Trying to solve DE ...
specfun-info: DSolve computes
         E^(2*x*z - z^2)*C[1]
specfun-info: expression rearranged:
         E^((2*x - z)*z)*C[1]
specfun-info: Calculation of initial values...
specfun-info: C[1] == 1

          (2 x - z) z
Out[35]= E

In[36]:= nospecfunprint
```

Note that `nospecfunprint` turns off the verbatim mode.

You saw that `Convert` uses the builtin Mathematica procedure `DSolve` to solve the differential equation corresponding to the given series coefficients. This can be rather time consuming or might be without success. If you know the generating function in advance, then you don't have to solve a differential equation.

The following two statements prove, for example, that

$$\sum_{k=0}^{\infty} P_n(x)\, z^n = \frac{1}{\sqrt{1 - 2xz + z^2}}$$

modulo two initial values without using `DSolve`. Here the procedure `SimpleRE[f,x,a[k]]` finds the recurrence equation which is valid for the series coefficients of `f` (see e.g. [12])

```
In[37]:= SimpleRE[1/Sqrt[1-2*x*z+z^2],z,P[n]]

Out[37]= -((1 + n) P[n]) + (3 + 2 n) x P[1 + n] -

>        (2 + n) P[2 + n] == 0

In[38]:= HolonomicRE[LegendreP[n,x],P[n]]

Out[38]= (1 + n) P[n] - (3 + 2 n) x P[1 + n] +

>        (2 + n) P[2 + n] == 0
```

## 5  Application: The $Z$-Transform

The $Z$-Transform

$$\mathcal{Z}\{a_k\} = f(z) = \sum_{n=0}^{\infty} a_k z^{-k}$$

of a sequence $\{a_k\}$ is the discrete analogue of the Laplace Transform. This series converges in the region outside the circle $|z| = |z_0| = \limsup_{k \to \infty} \sqrt[k]{|a_k|}$.

The procedure `ZTransform[a,k,z]` computes the $Z$-Transform of the sequence $a_k$, whereas the procedure `InverseZTransform[f,z]` gives the inverse $Z$-Transform of the functions $f(z)$, i.e., it computes $a_k$. Internally these procedures use `Convert` and `PowerSeries`, respectively. As an application we give the examples

```
In[39]:= fun=ZTransform[1/k!,k,z]

              1/z
Out[39]= E

In[40]:= InverseZTransform[fun,z]

              1 k
             (-)
              z
Out[40]= sum[----, {k, 0, Infinity}]
              k!

In[41]:= fun=ZTransform[1/(2*k+1)!,k,z]
```

```
                              Sqrt[1/z]
                    -1           E
Out[41]= -------------------- + ----------
            Sqrt[1/z]      1           1
         2 E          Sqrt[-]   2 Sqrt[-]
                           z           z

In[42]:= InverseZTransform[fun,z]

              1  k
             (-)
              z
Out[42]= sum[----------, {k, 0, Infinity}]
             (1 + 2 k)!

In[43]:= InverseZTransform[z/(1+z+z^2),z]

                1       k   1 1 + k
            (-(-------))   (-)
                   -2        z
                1 + z
Out[43]= sum[----------------------, {k, 0, Infinity}]
                        -2
                    1 + z
```

## 6  Application: Feynman Diagrams

In [5], Fleischer and Tarasov gave the representation[2]

$$V(\alpha,\beta,\gamma) = (-1)^{\alpha+\beta+\gamma} \cdot \frac{\Gamma(\alpha+\beta+\gamma-d/2)\Gamma(d/2-\gamma)\Gamma(\alpha+\gamma-d/2)\Gamma(\beta+\gamma-d/2)}{\Gamma(\alpha)\Gamma(\beta)\Gamma(d/2)\Gamma(\alpha+\beta+2\gamma-d)(m^2)^{\alpha+\beta+\gamma-d}}$$

$$\cdot {}_2F_1\left(\begin{array}{c}\alpha+\beta+\gamma-d,\ \alpha+\gamma-d/2\\ \alpha+\beta+2\gamma-d\end{array}\bigg| z\right)$$

for the calculation of a certain Feynman diagram. The function $V(\alpha,\beta,\gamma)$ can be easily computed for $\alpha,\beta,\gamma \in \{0,1\}$, hence one would like to attribute the calculation for nonnegative integer arguments $\alpha,\beta,\gamma$ to these cases. The simplest way is to compute the function by a pure recurrence equation that increases only one of its arguments, e.g. $\beta$, and leaves the other arguments constant. By Zeilberger's algorithm, the package gives the pure recurrence equation

```
In[44]:= HolonomicRE[(-1)^(alpha+beta+gamma)*
         Gamma[alpha+beta+gamma-d/2]*Gamma[d/2-gamma]*
         Gamma[alpha+gamma-d/2]*Gamma[beta+gamma-d/2]/
         (Gamma[alpha]*Gamma[beta]*Gamma[d/2]*
         Gamma[alpha+beta+2*gamma-d]*(m^2)^
         (alpha+beta+gamma-d))*
         Hypergeometric2F1[alpha+beta+gamma-d,
         alpha+gamma-d/2,alpha+beta+2*gamma-d,z],V[beta]]
```

---

[2] I would like to thank Jochem Fleischer who informed me about a typographical error in formula (31) of [5].

```
Out[44]= (2 beta - d + 2 gamma)

>       (2 alpha + 2 beta - d + 2 gamma)

>        (2 + 2 alpha + 2 beta - d + 2 gamma) V[beta] +
                                                       2
>     2 beta (2 + 2 alpha + 2 beta - d + 2 gamma) m

>       (2 alpha + 2 beta - 2 d + 4 gamma + 2 z +

>          2 beta z - d z) V[1 + beta] +

>      8 beta (1 + beta) (1 + alpha + beta - d + gamma)
           4
>       m   z V[2 + beta] == 0
```

with respect to $\beta$, and similar recurrence equations can be computed w.r.t. $\alpha$ and $\gamma$.

# 7 Application: Polynomial Approximations of Solutions of Differential Equations

Series solutions of holonomic differential equations satisfy holonomic recurrence equations from which their coefficients efficiently can be calculated by iteration. This is the asymptotically fastest known algorithm for the given purpose. It can be invoked by the Mathematica procedure `SeriesSolution[DE,y[x],initialvalues,n]`. Here `DE` is a given differential equation in terms of `y[x]`, `initialvalues` is a list of initial values, and `n` is the order of the proposed approximation. For the Airy function we have for example

```
In[45]:= Timing[SeriesSolution[y''[x]-x*y[x]==0,y[x],
         {AiryAi[0],AiryAi'[0]},10]]

                               x
Out[45]= {0.3 Second, -(--------------) -
                          1/3       1
                         3     Gamma[-]
                                    3

            4                7
           x                x
>      ---------------- - ----------------- -
           1/3       1        1/3       1
       12 3    Gamma[-]   504 3    Gamma[-]
                    3                   3

             10
            x                     1
>      ------------------ + ------------- +
            1/3       1       2/3       2
       45360 3    Gamma[-]   3     Gamma[-]
                     3                   3
```

$$
> \quad \frac{x^3}{6 \cdot 3^{2/3} \, \Gamma[\tfrac{2}{3}]^2} + \frac{x^6}{180 \cdot 3^{2/3} \, \Gamma[\tfrac{2}{3}]^2} +
$$

$$
> \quad \frac{x^9}{12960 \cdot 3^{2/3} \, \Gamma[\tfrac{2}{3}]^2} \Big\}
$$

whereas the buitin `Series` command gives

```
In[46]:= Timing[Series[AiryAi[x],{x,0,10}]]
```

$$
\text{Out[46]} = \Big\{ 3.5 \text{ Second}, \ \frac{1}{3^{2/3} \, \Gamma[\tfrac{2}{3}]} - \frac{x}{3^{1/3} \, \Gamma[\tfrac{1}{3}]} +
$$

$$
> \quad \frac{x^3}{6 \cdot 3^{2/3} \, \Gamma[\tfrac{2}{3}]^2} - \frac{x^4}{12 \cdot 3^{1/3} \, \Gamma[\tfrac{1}{3}]} +
$$

$$
> \quad \frac{x^6}{180 \cdot 3^{2/3} \, \Gamma[\tfrac{2}{3}]^2} - \frac{x^7}{504 \cdot 3^{1/3} \, \Gamma[\tfrac{1}{3}]} +
$$

$$
> \quad \frac{x^9}{12960 \cdot 3^{2/3} \, \Gamma[\tfrac{2}{3}]^2} - \frac{x^{10}}{45360 \cdot 3^{1/3} \, \Gamma[\tfrac{1}{3}]} +
$$

$$
> \quad O[x]^{11} \ \Big\}
$$

Hence even for modest order the method can be rather fast. If the order is large, then the speedup is even more impressive as can be seen from the computations below. Here `Taylor[f,x,x0,n]` gives a Taylor approximation of order $n$ for $f(x)$ at the point of development $x_0$, by generating a holonomic differential equation for $f(x)$, and using `SeriesSolution`, if applicable.

```
In[47]:= Timing[Series[AiryAi[x],{x,0,100}];]

Out[47]= {33.1833 Second, Null}

In[48]:= Timing[Taylor[AiryAi[x],{x,0,100}];]

Out[48]= {2.03333 Second, Null}

In[49]:= Timing[Series[Sin[x]*Exp[x],{x,0,100}];]

Out[49]= {5.88333 Second, Null}

In[50]:= Timing[Taylor[Sin[x]*Exp[x],{x,0,100}];]

Out[50]= {1.11667 Second, Null}

In[51]:= Timing[Series[Sin[x]*Exp[x],{x,0,1000}];]

Out[51]= {8517.53 Second, Null}

In[52]:= Timing[Taylor[Sin[x]*Exp[x],{x,0,1000}];]

Out[52]= {7.63333 Second, Null}

In[53]:= ps=PS[Sin[x]*Exp[x],x]

              k/2   k       k Pi
             2     x   Sin[----]
                             4
Out[53]= sum[-----------------, {k, 0, Infinity}]
                    k!

In[54]:= Timing[ps/.{sum->Sum,Infinity->1000};]

Out[54]= {7.21667 Second, Null}
```

Note that even the closed form solution given in line 43 does not provide a faster way to compute the series: The evaluation of each single coefficient takes about the same time than the iterative computation of the coefficients used by `Taylor`.

## 8  Advanced Applications

In this section, we give some advanced applications. Without explicit mentioning, all computations give proofs modulo initial values.

The first one is *Dougall's identity*

$${}_7F_6\left(\begin{matrix} a, 1+\frac{a}{2}, b, c, d, 1+2a-b-c-d+n, -n \\ \frac{a}{2}, 1+a-b, 1+a-c, 1+a-d, b+c+d-a-n, 1+a+n \end{matrix} \middle| 1\right)$$
$$= \frac{(1+a)_n \, (1+a-b-c)_n \, (1+a-b-d)_n \, (1+a-c-d)_n}{(1+a-b)_n (1+a-c)_n \, (1+a-d)_n \, (1+a-b-c-d)_n} \tag{2}$$

which is proved by

```
In[55]:= HolonomicRE[sum[HyperTerm[
         {a,1+a/2,b,c,d,1+2*a-b-c-d+n,-n},
         {a/2,1+a-b,1+a-c,1+a-d,b+c+d-a-n,1+a+n},1,k],
         {k,-Infinity,Infinity}],S[n]]

Out[55]= (1 + a + n) (1 + a - b - c + n)

>        (1 + a - b - d + n) (1 + a - c - d + n) S[n] -

>        (1 + a - b + n) (1 + a - c + n) (1 + a - d + n)

>        (1 + a - b - c - d + n) S[1 + n] == 0
```

Note that from this computation the parameters of the right hand hypergeometric term in (2) can be directly read off. Hence Zeilberger's algorithm has *discovered* the right hand side of (2) from its left hand input. Here the Mathematica procedure `HyperTerm[plist,qlist,x,k]` denotes the $k$th summand of the series `HypergeometricPFQ[plist,qlist,x]`.

Similarly *Clausen's identity*

$$_4F_3\left(\begin{array}{c} a,b,1/2-a-b-n,-n \\ 1/2+a+b,1-a-n,1-b-n \end{array}\bigg| 1\right)$$

$$= \frac{(2a)_n\,(a+b)_n\,(2b)_n}{(2a+2b)_n\,(a)_n\,(b)_n}$$

is proved by the computation

```
In[56]:= HolonomicRE[HypergeometricPFQ[
         {a,b,1/2-a-b-n,-n},
         {1/2+a+b,1-a-n,1-b-n},1],S[n]]

Out[56]= -((2 a + n) (a + b + n) (2 b + n) S[n]) +

>        (a + n) (b + n) (2 a + 2 b + n) S[1 + n] == 0
```

The only case when the square of a $_2F_1$ function is a $_3F_2$ function is given by *Clausen's formula*

$$_2F_1\left(\begin{array}{c} a,b \\ a+b+1/2 \end{array}\bigg| x\right)^2 = {_3F_2}\left(\begin{array}{c} 2a,2b,a+b \\ a+b+1/2,2a+2b \end{array}\bigg| x\right). \qquad (3)$$

Rather than proving this in the usual way by showing that both sides satisfy the same holonomic differential equation, we can use Zeilberger's algorihm to *determine* the right hand side of (3), given the left hand product, by using the Cauchy product. This is done by the computation

```
In[57]:= HolonomicRE[sum[
         HyperTerm[{a,b},{a+b+1/2},1,j]*
         HyperTerm[{a,b},{a+b+1/2},1,k-j],
         {j,0,k}],a[k]]
```

```
Out[57]= -2 (2 a + k) (a + b + k) (2 b + k) a[k] +

>       (1 + k) (2 a + 2 b + k) (1 + 2 a + 2 b + 2 k)

>       a[1 + k] == 0
```

Observe that one can directly read off the parameters of the hypergeometric function on the right hand side of (3).

The following identity of Askey and Gasper

$$\frac{(\alpha+2)_n}{n!} \,_3F_2\left(\begin{matrix} -n, n+\alpha+2, \frac{\alpha+1}{2} \\ \alpha+1, \frac{\alpha+3}{2} \end{matrix} \,\bigg|\, x\right) =$$

$$\sum_{j=0}^{[n/2]} \frac{\left(\frac{1}{2}\right)_j \left(\frac{\alpha}{2}+1\right)_{n-j} \left(\frac{\alpha+3}{2}\right)_{n-2j} (\alpha+1)_{n-2j}}{j! \left(\frac{\alpha+3}{2}\right)_{n-j} \left(\frac{\alpha+1}{2}\right)_{n-2j} (n-2j)!} \quad (4)$$

$$\cdot \,_3F_2\left(\begin{matrix} 2j-n, n-2j+\alpha+1, \frac{\alpha+1}{2} \\ \alpha+1, \frac{\alpha+2}{2} \end{matrix} \,\bigg|\, x\right)$$

(see e.g. [7]) was an essential tool in the proof of the *Bieberbach conjecture* by de Branges [3]. The Askey-Gasper identity is proved by the computation

```
In[58]:= HolonomicRE[sum[Pochhammer[1/2,j]*
        Pochhammer[alpha/2+1,n-j]*Pochhammer[(alpha+3)/2,n-2*j]*
        Pochhammer[alpha+1,n-2*j]/Pochhammer[(alpha+3)/2,n-j]/
        Pochhammer[(alpha+1)/2,n-2*j]/(n-2*j)!/j!*
        HyperTerm[{2*j-n,n-2*j+alpha+1,(alpha+1)/2},
        {alpha+1,(alpha+2)/2},1,k],{j,-Infinity,Infinity}],a[k]]

Out[58]= (1 + alpha + 2 k) (k - n) (2 + alpha + k + n)

>       a[k] - (1 + k) (1 + alpha + k)

>       (3 + alpha + 2 k) a[1 + k] == 0
```

Note that here the input is much more complicated than the output is! The computational trick is that we interchanged the order of summation of the double sum.

By an application of Clausen's formula from the Askey-Gasper identity it follows that for $\alpha > -2$ the left hand function in (4) is nonnegative. This was the result which de Branges needed for $\alpha = 0, 2, \ldots$ for his proof of the Bieberbach conjecture.

In recent work the *parameter derivatives* for the *Jacobi*, *Gegenbauer* and *Laguerre polynomials* have been determined ([6], [14]). For the Laguerre polynomials $L_n^{(d)}(x)$ one gets for example the simple formula

$$\frac{\partial}{\partial d} L_n^{(d)}(x) = \sum_{k=0}^{n-1} \frac{1}{n-k} L_k^{(d)}(x) \,.$$

This result can be obtained, e.g., by taking $\mu \to 0$ in the connection relation

$$L_n^{(\alpha+\mu)}(x) = \sum_{k=0}^{n} \frac{(\mu)_{n-k}}{(n-k)!} L_k^{(\alpha)}(x) \,.$$

The latter statement can be discovered by

```
In[59]:= HolonomicRE[sum[Pochhammer[mu,n-k]/(n-k)!*
         (-1)^j/j!*Binomial[k+alpha,k-j]*x^j,
         {k,-Infinity,Infinity}],a[j]]

Out[59]= (-j + n) x a[j] + (1 + j)
```
>       (1 + alpha + j + mu) a[1 + j] == 0

How can we obtain the exponential generating function

$$\sum_{n=0}^{\infty} \frac{1}{n!} P_n(x) z^n = e^{xz} J_0\left(z\sqrt{1-x^2}\right) \tag{5}$$

of the Legendre polynomials, $J_n(x)$ denoting the *Bessel functions*? Note that $J_0(x)$ is a ${}_0F_1$ hypergeometric function

```
In[60]:= FunctionToHypergeometric[BesselJ[0,x],x]

                              2
                            -x
Out[60]= hypergeometricPFQ[{}, {1}, ---]
                             4
```

To obtain (5), we use the representation

$$P_n(x) = x^n \, {}_2F_1\left(\begin{array}{c}-n/2, (1-n)/2\\ 1\end{array}\middle| 1 - \frac{1}{x^2}\right)$$

for the Legendre polynomials, and change the order of summation:

```
In[61]:= HolonomicRE[sum[x^n*
         HyperTerm[{-n/2,(1-n)/2},{1},1-1/x^2,k]/n!*z^n,
         {n,-Infinity,Infinity}],a[k]]

                 2           2
Out[61]= (1 - x) x  (1 + x) z  a[k] +

                 2  2
>       4 (1 + k)  x  a[1 + k] == 0
```

Now we can read off the hypergeometric parameters of the Bessel function, and using the initial value

$$\sum_{n=0}^{\infty} \frac{x^n}{n!} z^n = e^{xz}$$

we obtain (5).

# 9 Efficient Computation of Orthogonal Polynomials

In recent studies we discovered that none of the popular computer algebra systems Axiom, Macsyma, Maple, Mathematica, Mupad and REDUCE had implemented the most efficient algorithms to compute specific orthogonal polynomials, neither for rational exact nor for numerical computations, although Mathematica came rather near [16].

Hence, in the package `SpecialFunctions`, we have included efficient implementations for this purpose. If you compute the 1,000th Chebyshev polynomial by `Timing[ChebyshevT[1000,x];]` before and after loading the package, you will see how the computation is accelerated using the package. The larger $n$ is, the larger the advantage. (If you have a fast computer and enough memory, you should try $n = 5,000$ or $n = 10,000$.)

For $n = 0, 1, 2, \ldots$ the package supports the computation of the *classical orthogonal polynomials* `JacobiP[n,`$\alpha$`,`$\beta$`,x]`, `GegenbauerC[n,`$\alpha$`,x]`, `LaguerreL[n,`$\alpha$`,x]` and `HermiteH[n,x]`, as well as the computation of the *discrete orthogonal polynomials* (see [17]) `Hahn[n,N,`$\alpha$`,`$\beta$`,x]`, `DiscreteChebyshev[n,N,x]`, `Meixner[n,`$\gamma$`,`$\mu$`,x]`, `Krawtchouk[n,N,p,x]`, `DiscreteLaguerre[n,`$\rho$`,`$\alpha$`,x]` and `Charlier[n,`$\mu$`,x]`. Note that these as well as more special functions can also be used in connection with `HolonomicDE`, `HolonomicRE`, etc.

With `?SpecialFunctions`*`*` you get a list of all procedures exported by the package, and by `?PowerSeries`, e.g., you get a description of the `PowerSeries` procedure, together with an example. Moreover the contents of the integrated package `fastZeil` by Paule/Schorn can be listed by `?fastZeil`*`*`.

Some examples for the computation of specific orthogonal polynomials:

```
In[62]:= Charlier[5,mu,x]

              5 x     10 (1 - x) x
Out[62]= 1 - --- - ------------- -
             mu             2
                          mu

      10 (1 - x) (2 - x) x
  >   -------------------- -
              3
             mu

      5 (1 - x) (2 - x) (3 - x) x
  >   --------------------------- -
                 4
                mu

      (1 - x) (2 - x) (3 - x) (4 - x) x
  >   ---------------------------------
                      5
                     mu

In[63]:= ChebyshevT[1000,1/4]
```

```
Out[63]= 4633888252620729527555316242958909548932935 43\
>       931057925304868803092234816485137547165 3207603\
>       480107997324539642132129419798682743029733063451\
>       024106656353217059391151047974341312683 4882269\
>       510117040295963467578579678911640493759116734331\
>       488874299433540209147382435657218918439668245092\
>       8217512771529377 /
>       21430172143725346418968500981200036211 22809623411\
>       067214887500776740702102249872244986396757631391\
>       716255189345835106293650374290571384628087196915\
>       514939714960786913554964846197084214921012474228\
>       375590836430609294996716388253479753511833108789\
>       215412582914239295537308433532085966330524877367\
>       4411336138752

In[64]:= N[ChebyshevT[10^10,N[1/4,100]],100]

Out[64]= -0.1615907100588309730645451317842686501 11183\
>       5308168473991177711097513735275666631253 3037757
```

Note that for plotting purposes, Mathematica's original implementation is advantageous. Using the Charlier polynomials, we finally give some examples of the symbolic use of the discrete orthogonal families:

```
In[65]:= HolonomicRE[Charlier[n,mu,x],C[x]]

Out[65]= (-1 - x) C[x] + (1 + mu - n + x) C[1 + x] -
>       mu C[2 + x] == 0

In[66]:= HolonomicRE[Charlier[n,mu,x],C[n]]

Out[66]= (1 + n) C[n] + (-1 - mu - n + x) C[1 + n] +
>       mu C[2 + n] == 0

In[67]:= HolonomicRE[Binomial[mu,n]*Charlier[n,mu,x],a[n]]

Out[67]= (-mu + n) (1 - mu + n) a[n] +
```

```
>         (1 - mu + n) (1 + mu + n - x) a[1 + n] +

>       mu (2 + n) a[2 + n] == 0

In[68]:= HolonomicRE[Charlier[n,mu,2x],a[x]]

Out[68]= 2 (-3 - mu + n - 2 x) (1 + x) (1 + 2 x)

                                 2       3
>         a[x] + (6 + 2 mu + mu  + mu  - 11 n - 7 mu n -

              2         2            2       3
>         3 mu  n + 6 n  + 3 mu n  - n  + 22 x +

                   2
>         6 mu x + 2 mu  x - 24 n x - 8 mu n x +

             2        2         2        2        3
>         6 n  x + 24 x  + 4 mu x  - 12 n x  + 8 x  )

                     2
>         a[1 + x] + mu  (-1 - mu + n - 2 x) a[2 + x] == 0
```